\documentclass[11pt]{article}
\usepackage{amssymb}
\usepackage{amsmath}
\usepackage{times}\hfuzz=10pt 
   \csname @addtoreset\endcsname{equation}{section}
\newtheorem{theorem}{Theorem}%[section]
\newtheorem{lemma}{Lemma}%[section]
\newtheorem{corollary}{Corollary}%[section]
\newtheorem{definition}{Definition}%[section]
\def\e{\varepsilon}

\def\defi{\stackrel{{\scriptscriptstyle \Delta}}{=}}

\def\a{\alpha}
\def\d{\delta}
\def\o{\omega}
\def\O{\Omega}

\def\w{\widehat}
\def\Ind{{\mathbb{I}}}

\def\sign{{\rm  sign\,}}
\def\mes{{\rm mes\,}}

\def\Re{{\rm Re\,}}

\def\R{{\bf R}}

\def\Z{{\cal Z}}

\def\g{\gamma}
\def\C{{\bf C}}

\def\ww{\widetilde}
\def\X{{\cal X}}

\def\oo{\bar}

\def\G{\Gamma}

\def\T{{\mathbb{T}}}

\newcommand{\be}{\begin{equation}}
\newcommand{\ee}{\end{equation}}
\newcommand{\bd}{\begin{displaymath}}
\newcommand{\ed}{\end{displaymath}}
\newcommand{\ba}{\begin{array}{ll}}
\newcommand{\ea}{\end{array}}
\newcommand{\baa}{\begin{eqnarray}}
\newcommand{\eaa}{\end{eqnarray}}
\newcommand{\baaa}{\begin{eqnarray*}}
\newcommand{\eaaa}{\end{eqnarray*}}
\font\sm=cmr10
\def\oo{\bar}

\def\a{\alpha}

\def\K{{\cal K}}
\def\Ko{K}%{K_0}
\def\ko{k}%{K_0}
\def\Yo{Y}%{K_0}
\def\yo{y}%{K_0}
\def\defi{=}
\date{Submitted: December 1,2011. Revised: March 14, 2012 }
\date{1st version: December 1,2011. Revised version: March 15, 2012 }%arcChiV
\title{On predictors for  band-limited and high-frequency time series\footnote { {\em Signal Processing} {\bf 92}, iss.
10, pp.2571-2575.   In {\em Fast communications} section.
DOI:10.1016/j.sigpro.2012.04.006}}
\author{
Nikolai Dokuchaev \\ {\sm  Department of Mathematics and Statistics,
Curtin University},\\ {\sm GPO Box U1987, Perth, Western Australia,
6845}\\ {\sm email N.Dokuchaev@curtin.edu.au. Tel.: 61 8 92663144.}}
 
\begin{document}
 \vspace{-0.5cm}      \maketitle
\begin{abstract} Pathwise predictability  and predictors for discrete time processes are studied in deterministic setting.
It is suggested to approximate convolution sums over future times by
convolution sums over past time. It is shown that all band-limited
processes are predictable in this sense, as well as high-frequency
processes with zero energy at low frequencies. In addition, a
process of mixed type still can be predicted if an ideal low-pass
filter exists for this process.
\\{
PACS 2008 numbers: 02.30.Mv, %    Approximations and expansions
02.30.Nw,  %  Fourier analysis
02.30.Yy, %    Control theory
07.05.Mh,  %  Neural networks, fuzzy logic, artificial intelligence
07.05.Kf %Data analysis: algorithms and implementation; data management
}\end{abstract}
\section{Introduction}  We study pathwise predictability of
discrete time processes in deterministic setting. It is well known
that certain restrictions on frequency distribution can ensure
additional opportunities for prediction and interpolation of the
processes. The classical result is Nyquist-Shannon-Kotelnikov
interpolation theorem for the continuous time band-limited
processes. It is also known that optimal prediction error for
stationary Gaussian processes is zero for the case of degenerate
spectral density. The related results can be found in Wainstein and
Zubakov (1962), Knab (1981), Papoulis (1985), Marvasti (1986),
Vaidyanathan (1987), Lyman {\it et al} (2000, 2001), Dokuchaev
(2008,2010).
\par
The present paper extends on discrete time setting the approach
suggested for continuous time processes in Dokuchaev (2008). We
study a special kind of predictors such that convolution sums over
future are approximated by convolution sums over past times
representing historical observations. We found some cases when this
approximation can be made uniformly over a wide class of input
processes, including all band-limited processes and high-frequency
processes.
 For the
processes of mixed type, we found that the similar predictability
can be achieved when the model allows a low pass filter  that acts
as an ideal low-pass filter for this process. These results can be
a useful addition to the existing theory of band-limited
processes. The novelty is that we consider predictability of both
high frequent and band-limited processes in a weak sense uniformly
over classes of input processes. In addition, we suggest a new
type of predictor. Its kernel is given explicitly in the frequency
domain.
\section{Definitions}
 Let  $D\defi\{z\in\C:
|z|\le 1\}$, $D^c\defi \C\setminus D$, $\T\defi\{z\in\C:\ |z|=1\}$,
\par
We denote by $\ell_r$ the set of all sequences
$x=\{x(t)\}_{t=-\infty}^{\infty}\subset\C$ such that
$\|x\|_{\ell_r}=\left(\sum_{t=-\infty}^{\infty}|x(t)|^r\right)^{1/r}<+\infty$
for $r\in[1,\infty)$,  $\|x\|_{\ell_\infty}=\sup_t|x(t)|<+\infty$
for $r=+\infty$. \par Let $\ell_r^+$ be the set of all sequences
$x\in\ell_r$ such that $x(t)=0$ for $t=-1,-2,-3,...$.

\par
For complex valued sequences $x\in \ell_1$ or $x\in \ell_2$, we
denote by $X=\Z x$  the Z-transform  \baaa
X(z)=\sum_{t=-\infty}^{\infty}x(t)z^{-t},\quad z\in\C. \eaaa
Respectively, the inverse $x=\Z^{-1}X$ is defined as \baaa
x(t)=\frac{1}{2\pi}\int_{-\pi}^\pi X\left(e^{i\o}\right) e^{i\o
t}d\o, \quad t=0,\pm 1,\pm 2,....\eaaa

 If $x\in \ell_2$, then $X|_\T$ is defined as an element of
$L_2(\T)$.
\par
Let  $H^r$ be the Hardy space of functions that are holomorphic on
$D^c$ including the point at infinity   (see, e.g., Duren (1970)).
Note that Z-transform defines bijection between the sequences from
$\ell_2^+$ and the restrictions (i.e. traces) of the functions from
$H^2$  on $\T$.
\begin{definition}
Let $\K$ be the class of all functions $k\in\ell_{\infty}$ such that
$k (t)=0$ for $t>0$ and $K=\Z k$ is \be K(z)=\frac{d(z)}{\d(z)},
\label{kda}\ee where $d(\cdot)$ and $\d(\cdot)$ are polynomials such
that ${\rm deg\,} d < {\rm deg\,} \d$, and if $\d(z)=0$ for $z\in\C$
then $|z|>1$.
\end{definition}
The class includes all kernels $k$ representing the anti-causal
linear constant-coefficient difference equations.
\begin{definition}
Let $\w\K$ be the class of functions $\w k:\ell_{\infty}^+$ such
that the function $\w K(\cdot)=\Z\w k$ belongs to $H^\infty\cap
H^2$.
\end{definition}
\par It follows from the definitions that if
 $\w k\in \w\K$ then $\w k (t)=0$ for $t<0$.
\par
 We are going to study linear predictors in the form $\w y(t)=\sum_{s=-\infty}^t\w
k(t-s)x(s)$ for the processes  $y(t)=\sum_{s=t}^{+\infty}
k(t-s)x(s)$, where $k\in\K$ and $\w k\in\w K$. The predictors use
historical values of currently observable process $x(\cdot)$.
\par
\begin{definition}\label{def1}
Let  $\X=\{x(\cdot)\}$ be a class of sequences from $\ell_{\infty}$,
let $r\in[1,+\infty]$, and let $\ww\K\subset\K$ be a class of
sequences.
\begin{itemize}
\item[(i)]
 We say that the class $\X$ is  $\ell_r$-predictable in
the weak sense with respect to the class $\ww\K$ if, for any
$k(\cdot)\in\ww\K$, there exists a sequence $\{\w
k_m(\cdot)\}_{m=1}^{+\infty}=\{\w
k_m(\cdot,\X,k)\}_{m=1}^{+\infty}\subset \w\K$ such that $$ \|y-\w
y_m\|_{\ell_r}\to 0\quad \hbox{as}\quad m\to+\infty\quad\forall
x\in\X,
$$ where \baaa y(t)\defi \sum_{s=t}^{+\infty}k(t-s)x(s),\qquad \w
y_m(t)\defi \sum^t_{s=-\infty}\w k_m(t-s)x(s).\label{predict} \eaaa
\item[(ii)] Let the set $\Z(\X)\defi \{X(e^{i\o})=\Z x|_{\T},\quad x\in \X\}$  be provided with a norm $\|\cdot\|$.
 We say that the class $\X$ is  $\ell_r$-predictable in
the weak sense  with respect to the class $\ww\K$ uniformly with
respect to the norm $\|\cdot\|$, if, for any $k(\cdot)\in\ww\K$ and
$\e>0$, there exists $\w k(\cdot) =\w k (\cdot,\X,k,\|\cdot\|,\e)\in
\w\K$ such that $$ \|y- \w y\,\|_{\ell_r}\le \e\|X\|\quad \forall
x\in\X,\quad X=\Z x. $$ Here $y(\cdot)$ is the same as above, $\w
y(t)\defi \sum^t_{s=-\infty}\w k(t-s)x(s).$
\end{itemize}
\end{definition}
 We call functions $\w k(\cdot)$ in Definition \ref{def1} predictors or
predicting kernels.
\section{The main result}
\par
Let $\O\in(0,\pi)$ be given, and let \baaa \X_L\defi\{x(\cdot)\in
\ell_2:\
X\left(e^{i\o}\right)=0\quad\hbox{if}\quad|\o|>\O,\quad X=\Z x\},\\
\X_H\defi\{x(\cdot)\in \ell_2:\
X\left(e^{i\o}\right)=0\quad\hbox{if}\quad|\o|<\O,\quad X=\Z x\}.
\eaaa
 In particular, $\X_L$ is a
class of band-limited processes, and $\X_H$ is a class of
high-frequency processes.
\subsection{Predictability of band-limited and high-frequency
processes from $L_2$} Let $\K_0$ be the class of all functions
$k\in\ell_{\infty}$ such that  $k (t)=0$ for $t>0$ and that $K=\Z k$
can be represented as \be K(z)=\frac{z+b}{z+a}, \label{kda1}\ee for
some real $a\in (-\infty,-1)\cup(1,+\infty)$ and $b\in\R$.
\begin{theorem}\label{ThM}\begin{itemize}\item[(i)] The classes $\X_L$ and
$\X_H$ are $\ell_2$-predictable in the weak sense  with respect to
the class $\K_0$.
\item[(ii)] The
classes $\X_L$ and $\X_H$  are $\ell_\infty$-predictable in the weak
sense  with respect to the class $\K_0$ uniformly with respect to
the norm $\|X(e^{i\o})\|_{L_2(-\pi,\pi)}$.
\item[(iii)] For any $q>2$, the classes $\X_L$ and $\X_H$  are
$\ell_2$-predictable  in the weak sense with respect to the class
$\K_0$ uniformly with respect to the norm
$\|X(e^{i\o})\|_{L_q(-\pi,\pi)}$.
\end{itemize}
\end{theorem}
\par  The question arises how to find the predicting kernels. In
the proof of Theorem \ref{ThM}, a possible choice of the kernels is
given explicitly via  Z-transforms.
\section{On a model with ideal low pass-pass filter}
\begin{corollary}\label{cor1}
Assume a model with  a process $x(\cdot)$ such that it is possible
to decompose it
 as $x(t)=x_L(t)+x_H(t)$, where   $x_L(\cdot)\in\X_L$ and
$x_H(\cdot)\in\X_H$. Then this observer would be able to predict
(approximately, in the sense of weak predictability with respect to
the class $\K_0$) the values of $y(t)=
\sum_{s=t}^{+\infty}k(t-s)x(s) $ for $k(\cdot)\in\K$ by predicting
the processes $y_L(t)=\sum_{s=t}^{+\infty}k(t-s)x_L(s) $ and
$y_H(t)=\sum_{s=t}^{+\infty}k(t-s)x_H(s) $ separately. More
precisely, the process $\w y(t)\defi \w y_L(t)+\w y_H(t)$ is the
prediction of $ y(t)$, where $y_L(t)=\sum^t_{-\infty}\w
k_L(t-s)x_L(s) $ and $y_H(t)=\sum^t_{-\infty}\w k_H(t-s)x_H(s) $,
and where $\w k_L(\cdot)$ and $\w k_H(\cdot)$ are predicting kernels
which existence for the processes $x_L(\cdot)$ and $x_H(\cdot)$  is
established above.
\end{corollary}
\par
Let $\chi_L\left(e^{i\o}\right)\defi \Ind_{\{|\o|\le \O\}}$ and
$\chi_H\left(e^{i\o}\right)\defi
1-\chi_L\left(e^{i\o}\right)=\Ind_{\{|\o|> \O\}}$, where $\o\in\R$;
$\Ind$ denote the indicator function.
\par
 The assumptions of Corollary \ref{cor1} mean
that there are a low-pass filter  and a high-pass filter with the
transfer functions $\chi_L$ and $\chi_H$ respectively, with
$x(\cdot)$ as the input, i.e., that the values $x_L(s)$ and $x_H(s)$
for $s\le t$ are available at time $t$, where \baaa x_L(\cdot)\defi
\Z^{-1}X_L,\quad X_L\left(e^{i\o}\right)\defi \chi_
L\left(e^{i\o}\right)X\left(e^{i\o}\right),\\ x_H(\cdot)\defi
\Z^{-1}X_H,\quad X_H\left(e^{i\o}\right)\defi
\chi_H\left(e^{i\o}\right)X\left(e^{i\o}\right),\eaaa and where
$X\defi \Z x$.
 It
follows that the predictability in the weak sense with respect to
the class $\K_0$ is possible for any process $x(\cdot)$ that can be
decomposed without error on a band limited process and a
high-frequency process, i.e., when there is a low-pass filters which
behave as an ideal filter for this process. (Since
$x_H(t)=x(t)-x_L(t)$, existence of the low pass filter implies
existence of the high pass filter). On the other hand, Corollary
\ref{cor1} implies that the existence of ideal low-pass filters is
impossible for general processes, since they cannot be predictable
in the sense of Definition \ref{def1}.
\par
Clearly, processes $x(\cdot)\in\X_L\cup \X_H$ are automatically
covered by Corollary \ref{cor1}, i.e., the existence of the filters
is not required for this case. For instance, we have immediately
that $x_L(\cdot)=x(\cdot)$ and $x_H(\cdot)\equiv 0$ for band-limited
processes.
\section{Proofs}
It suffices to present a set of predicting kernels $\w k$ with the
desired properties. We will use a version of the construction
introduced in Dokuchaev (2008) for continuous time setting.  This
construction is very straightforward and does not use the advanced
theory of $H^p$-spaces.
\par
Let $\K_1$ be the class of all functions $k\in\K_0$ such that $K=\Z
k$ can be represented as \baa K(z)=\frac{1}{z+a}, \label{kda3}\eaa
for some real $a\in (-\infty,-1)\cup(1,+\infty)$.

If $k\in\K_0$, then $K=\Z k$ can be represented as \baaa
K(z)=\frac{z+b}{z+a}= \frac{z+a+b-a}{z+a}=1+\frac{c}{z+a},
\label{kda2}\eaaa with $a\in (-\infty)\cup(1,+\infty)$, $b\in\R$,
and $c=b-a$. It follows that the process $y(t)$ for $k\in\K_0$ can
be represented as $y(t)=x(t)+c\sum_{s=t}^{+\infty}k_1(t-s)x(s)$,
where $k_1\in\K_1$. Therefore, it suffices to prove theorem for
$k\in\K_1$ only.

Let  $k(\cdot)\in\K_1$ and $K\left(e^{i\o}\right)=\Z k$ be defined
by (\ref{kda3}) for some for  $a\in (-\infty,-1)\cup(1,+\infty)$.

Let $G=(-\O,\O)$, and let  \baa \a=-\frac{1+a\cos(\O)}{a+\cos(\O)}.
\label{sa}\eaa Let us show that $\a=f(a)\in (-1,1)$. Clearly, the
function
$$f(a)=\frac{1+a\cos(\O)}{a+\cos(\O)}$$ is such that $f'(a)<0$ for all $a$ such that $|a|\ge 1$, $f(-1)=-1$, $f(1)=1$, and $f(\pm\infty)=\cos(\O)$.
These properties imply that $\a=f(a)\in (-1,1)$.
\par
Further, we have that $1+\a a+(a+\a)\cos(\O)=0, $ and \baa
&&\sign(a+\a)(1+\a a+(a+\a)\cos(\o))>0,\quad \o\in G,\nonumber\\
&&\sign(a+\a)(1+\a a+(a+\a)\cos(\o))<0,\quad \o\in(-\O,\O)\backslash
G. \label{ineqa}\eaa
\par
  %(26 September 2008)
Set \baa V(z)\defi
1-\exp\left(\g\sign(a+\a)\frac{z+a}{z+\a}\right),\quad
  \w \Ko (z)\defi V(z)\Ko (z),\quad \g\in\R.
  \label{K}\eaa

\begin{lemma}\label{lemmaV}
\begin{itemize}
\item[(i)] $V(z)\in H^{\infty}$ and $\w \Ko (z)\defi \Ko(z)V(z)\in
H^{\infty}\cap H^2$.
\item[(ii)]  If $\g<0$ and $\o\in[-\O,\O]$, then  $|V(e^{i\o})|\le 2$.
If $\g>0$ and $\o\in[-\pi,\pi]\backslash(-\O,\O)$, then
$|V(e^{i\o})|\le 2$.
\item[(iii)] If $\o\in(-\O,\O)$, then $V(e^{i\o})\to 1$ as
$\g\to -\infty$. If $\o\in [-\pi,\pi]\backslash[-\O,\O]$, then
$V(e^{i\o})\to 1$ as $\g\to +\infty$.
\item[(iv)]
For any $\e\in(0,\O)$, $V(e^{i\o})\to 1$ as $\g\to-\infty$ uniformly
in $\o\in[-\O+\e,\O-\e]$ as $\g\to-\infty$, and  $V(e^{i\o})\to 1$
as $\g\to +\infty$ uniformly in $\o\in
[-\pi,\pi]\backslash(-\O+\e,\O-\e)$. \index{
\item[(v)] The function $\w K\left(e^{i\o}\right)^{-1}$ is bounded on
$\o\in[-\pi,\pi]$.}
\end{itemize}
\end{lemma}
\par
{\it Proof of Lemma \ref{lemmaV}}. Clearly, $V\in H^{\infty}$, and
$(z+a)^{-1}V(z)\in H^2\cap H^{\infty}$, since the pole of
$(z+a)^{-1}$ is being compensated by multiplying with $V$. It
follows that $K(z)V(z)\in H^2\cap H^{\infty}$. Then statement (i)
follows.
\par Further, for $\o\in\R$,
 \baaa
\frac{e^{i\o}+a}{e^{i\o}+\a}=\frac{(e^{i\o}+a)(e^{-i\o}+\a)}{|e^{i\o}+\a|^2}=\frac{1+a\a+ae^{-i\o}+\a
e^{i\o}}{|e^{i\o}+\a|^2}.\eaaa Hence
 \baaa
\Re\frac{e^{i\o}+a}{e^{i\o}+\a}
=\frac{1+a\a+(a+\a)\cos(\o)}{|e^{i\o}+\a|^2}. \label{re} \eaaa
 Then statements (ii)-(iv)  follow from (\ref{ineqa}). \index{The proof of statement (v) is
 straightforward.}
 This completes the proof of Lemma \ref{lemmaV}. $\Box$
\par
{\it Proof of Theorem \ref{ThM}}.  For  $x(\cdot)\in \ell_2$, let
$X\defi \Z x$, $\ko=\Z^{-1}\Ko $, $\w \ko=\Z^{-1}\w \Ko$,   \baaa
\yo (t)\defi \sum_{s=t}^{\infty}\ko (t-s)x(s),\quad \w \yo (t)\defi
\sum^t_{s=-\infty}\w \ko   (t-s)x(s). \eaaa Let
 $\Yo \defi \Z \yo $, let $V$ and $\w K$ be as defined above, and let
 $\w Y\defi \w K X$.
\par
 Let us consider the cases of $\X_L$ and $\X_H$ simultaneously.
For the case of the class $\X_L$, consider $\g<0$ and assume that
$\g\to -\infty$. Set $\G\defi [-\O,\O]$ for this case. For the case
of the class $\X_H$, consider $\g>0$ and $\g\to +\infty$. Set
$\G\defi [-\pi,-\O]\cup[\O,+\pi]$ for this case.
\par
Let $x(\cdot)\in\X_L$ or $x(\cdot)\in\X_H$. In  both  cases, Lemma
\ref{lemmaV} gives that $|V\left(e^{i\o}\right)|\le 2$ for all
$\o\in \G$. If $\g\to -\infty$ or $\g\to +\infty$ respectively  for
$\X_L$ or $\X_H$ cases, then
  $V\left(e^{i\o}\right)\to 1$  for a.e. $\o\in \G$, i.e., for a.e. $\o$ such that
  $X\left(e^{i\o}\right)\neq 0$.
\par
Let us prove (i).  Since $K\left(e^{i\o}\right)\in
L_{\infty}(-\pi,\pi)$, $\w K\left(e^{i\o}\right)\in
L_{\infty}(-\pi,\pi)$, and $X\left(e^{i\o}\right)\in L_2(-\pi,\pi)$,
we have that
$Y\left(e^{i\o}\right)=K\left(e^{i\o}\right)X\left(e^{i\o}\right)\in
L_{2}(-\pi,\pi)$ and $\w Y\left(e^{i\o}\right)=\w
K\left(e^{i\o}\right)X\left(e^{i\o}\right)\in L_{2}(-\pi,\pi)$. By
Lemma \ref{lemmaV}, it follows that  \baa\w Y\left(e^{i\o}\right)\to
Y\left(e^{i\o}\right)\quad\hbox{for a.e.}\quad \o\in \R,
\label{YY}\eaa
 as $\g\to -\infty$
or $\g\to +\infty$ respectively for $\X_L$ or $\X_H$ cases. We have
that \baa&&|\w K\left(e^{i\o}\right)-K\left(e^{i\o}\right)|\le
|V\left(e^{i\o}\right)-1||K\left(e^{i\o}\right)|\le 2
|K_m\left(e^{i\o}\right)|,\quad \o\in \G,\label{d1}\\ &&|\w
Y\left(e^{i\o}\right)-Y\left(e^{i\o}\right)|\le
2|Y\left(e^{i\o}\right)|=2|K\left(e^{i\o}\right)||X\left(e^{i\o}\right)|,\quad
\o\in \G.\label{d2} \eaa
 By (\ref{YY}),(\ref{d2}), and by Lebesque Dominance Theorem, it follows that
\be\|\w
Y\left(e^{i\o}\right)-Y\left(e^{i\o}\right)\|_{L_2(-\pi,\pi)}\to
0,\quad\hbox{i.e.,}\quad\|\w y-y\|_{L_2(-\pi,\pi)}\to 0
\label{1s}\ee
 as $\g\to -\infty$
or $\g\to +\infty$ respectively for $\X_L$ or $\X_H$ cases, where
$\w y=\Z^{-1}\w Y$.
\par
 Let us prove (ii)-(iii).
 Take $d=1$ for (ii) and take $d=2$ for (iii). If
$X\left(e^{i\o}\right)\in L_\nu(-\pi,\pi)$ for $\nu>d$, then
H\"older inequality gives
 \be\|\w Y\left(e^{i\o}\right)-Y\left(e^{i\o}\right)\|_{L_d(-\pi,\pi)}\le
\|\w
K\left(e^{i\o}\right)-K\left(e^{i\o}\right)\|_{L_\mu(\G)}\|X\left(e^{i\o}\right)\|_{L_\nu(\G)},
\label{4s}\ee where $\mu$ is such that $1/\mu+1/\nu=1/d$.
 By (\ref{d1}) and by Lebesque Dominance Theorem again, it follows that
\be\label{2s}\|\w
K\left(e^{i\o}\right)-K\left(e^{i\o}\right)\|_{L_\mu(\G)}\to 0\quad
\forall \mu \in[1,+\infty), \ee  as $\g\to -\infty$ or $\g\to
+\infty$ respectively for $\X_L$ or $\X_H$ cases. Then, by
(\ref{4s})-(\ref{2s}), it follows  that the predicting kernels $\w
k(\cdot)=\w k(\cdot,\g)=\Z^{-1}\w K$ are such as required in
statements (ii)--(iii). This completes the proof of Theorem
\ref{ThM}. $\Box$

\par Corollary \ref{cor1} follows
immediately from Theorem \ref{ThM}.
\section{On the prediction error generated by a  high-frequency noise}
Let us estimate the prediction error  for the case when predictor
(\ref{K}) designed for a band-limited process is applied to a
process with a small high-frequency noise.

Let $\O\in(0,\pi)$  and $\nu\in [0,1)$  be given. Let us consider a
process $x(\cdot)\in\ell_{\infty}$ such that  $|X(i\o)|\le 1$ for
$\o\in G$ and
 $|X(i\o)|\le \nu$ for $\o\in[-\pi,\pi]\backslash G$ , where $X=\Z x$ and $G=(-\O,\O)$.

Assume that predictor  (\ref{K}) is constructed under the hypothesis
that $\nu=0$ (i.e, that $x(\cdot)$ is a band-limited processes from
$\X_L$), for some $a\in\R\backslash [-1,1]$.   For an arbitrarily
small $\e>0$, we can find $\g=\g(\e)$ such that if the hypothesis
that $\nu=0$ is correct, then \baa \|\w
y-y\|_{\ell_{\infty}}\le\frac{\e}{2\pi},\label{eps}\eaa  where
$y(\cdot)$ and $\w y(\cdot)$ are such as in Definition \ref{def1}.
\par
Let us estimate the prediction error for the case when $\nu>0$. We
have that \baaa \|\w y-y\|_{\ell_{\infty}}\le\frac{1}{2\pi}\|\w
Y\left(e^{i\o}\right)-Y\left(e^{i\o}\right)\|_{L_1(-\pi,\pi)}, \eaaa
where $Y=\Z y$ and $\w Y=\Z\w y$. Let $\O_1=\O-\e/4$ and
$G_1=(-\O_1,\O_1)$. By the assumptions on $X$, we have that \baaa
\|\w
Y\left(e^{i\o}\right)-Y\left(e^{i\o}\right)\|_{L_1(-\pi,\pi)}\le
I_1+I_2+\nu I_3, \eaaa where \baaa I_1=\kappa
\int_{G_1}e^{\g\psi(\o)}d\o,\quad I_2=\kappa\int_{G\backslash
G_1}e^{\g\psi(\o)}d\o,\quad I_3=\kappa\int_{(-\pi,\pi)\backslash
G}e^{\g\psi(\o)}d\o, \eaaa and where
$\kappa=\max_{\o}||K(e^{i\o})|$, \baaa
\psi(\o)=\sign(a+\a)\Re\frac{e^{i\o}+a}{e^{i\o}+\a}. \eaaa Note that
$\psi(\o)>0$ for $\o\in G$. Let $\psi_0=\min_{\o\in G_1} \psi(\o)$,
and let $\g=-\log(2\kappa/\e)/\psi_0$. Then $I_1\le \e/2$. Further,
$I_2\le \kappa\,\mes(G\backslash G_1)=\e/2$. Therefore, (\ref{eps})
holds if $\nu=0$.
\par
The value $I_1+I_2$ represents the forecast error when $\nu=0$; this
error can be done arbitrarily small with $\g$ selected as above when
$\e\to 0$.
\par
 Let us estimate $I_3$.
Clearly, $|\psi(\o)|\le 1+\left|\frac{a-\a}{e^{i\o}+\a}\right|\le
\mu$, where $\mu\defi 1+|a-\a|/(1-\a)$. Hence \baaa I_3\le
\kappa\int_{(-\pi,\pi)\backslash G}e^{\g
\mu}d\o=\kappa\int_{(-\pi,\pi)\backslash
G}e^{\frac{\log(2\kappa/\e)}{\psi_0}\mu}d\o=2\kappa(\pi-\O)e^{\frac{\log(2\kappa/\e)\mu}{\psi_0}}.\eaaa
Hence \baaa\nu I_3\le
2\kappa\nu(\pi-\O)\left(\frac{2\kappa}{\e}\right)^{\frac{\mu}{\psi_0}}.\eaaa
The value $\nu I_3$  represents the additional error caused by the
presence of unexpected high-frequency noise (when $\nu>0$). It can
be seen that if $\e\to 0$ than this error is increasing as a
polynomial of $\e^{-1}$ with the rate depending on $\a$ (defined by
$\O$ and $a$). If $\O\to \pi$ then $|\a|\to 1$ and $\mu\to +\infty$,
and, for a given $\e$,  the error is increasing exponentially in
$\mu$.
%\end{document}
\section{Concluding remarks}
\begin{itemize}
\item By (\ref{sa}),  $\a\to \pm 1$ as
$\O\to\pi$,  and the predictor suggested above loses its feasibility
as $\O\to\pi$. (In particular, $\|\w k\|_{\ell_{\infty}}\to
+\infty$).
\item If $k(\cdot)$ is a real valued function, then $\w k$ is also
real valued. It follows from the fact that $K\left(\oo
z\right)=\overline{K\left(z\right)}$, and, therefore,
$K\left(e^{-i\o}\right)=\overline{K\left(e^{i\o}\right)}$.
\item A similar approach  can be applied to the case when
 $X(z)$ vanishes on some connected set $I\subset\T$. In this case,
 the classes $\Ko_0$ and $\Ko_1$ have to be replaced by similar
 classes with complex  $a\in D^c$. For real valued kernels,
 it could be meaningful to include the functions $K$
 represented by the sums of two simple fractions, to ensure that the
 process
 $\Z^{-1}k$ is real (i.e, that $\overline{K\left(e^{i\o}\right)}=K\left(e^{-i\o}\right))$.
\item
 The predictors obtained above require the past values of $x(s)$
for all $s\in(-\infty,t]$. In practice, $\sum_{s=-\infty}^t\w
k(t-s)x(s) $ can be approximated by $\sum_{s=-M}^t\w k(t-s)x(s) $
for large enough $M>0$. In addition, the corresponding transfer
functions can be approximated by rational fraction polynomials.
\item
 The system for the suggested
predictors is stable, since the corresponding transfer functions
have poles in the domain $\{|z|<1\}$ only.  However, the suggested
predictors are not robust. For instance, if the predictor is
designed for the class $\X_L$ and it is applied for a process
$x(\cdot)\notin\X_L$ with small non-zero energy at the frequencies
outside $[-\O,\O]$, then the error generated by the presence of this
energy is increasing if $\g\to\infty$. \item The results of this
paper can be applied to discrete time stationary random Gaussian
processes.  In particular, assume that the spectral density of the
underlying process $x(t)$ vanishes outside the interval
$[-\O,\O]\subset (-\pi,\pi)$. It is known that the minimal (optimal)
predicting error is zero in this case. The sequence of the
predictors constructed above represents a sequence of suboptimal
predictors leading to vanishing prediction error.
\end{itemize}
\subsection*{Acknowledgments}  This work  was supported by ARC grant of Australia DP120100928 to the author.
In addition, the author thanks Prof. Augusto Ferrante for useful
discussion and advice on  spectral analysis of time series.

\end{document}